\documentclass[english]{article}
\usepackage[T1]{fontenc}
\usepackage[latin9]{inputenc}
\usepackage{float}
\usepackage{amsmath}
\usepackage{amssymb}
\usepackage{stackrel}
\usepackage{esint}
\usepackage{nameref}

\makeatletter

\usepackage{xr}
\usepackage{hyperref}
\usepackage[capitalise]{cleveref}

\usepackage{pgfplots}
\usepackage{pgfplotstable}
\usepgfplotslibrary{units}
\usepackage{booktabs}
\usepackage{authblk}
\usepackage[english]{babel}

\makeatother

\begin{document}

\title{An Elastic Energy Minimization Framework for Mean Surface Calculation}

\author[$1$]{Jozsef Molnar} 
\author[$1,2$]{Peter Horvath}

\affil[$1$]{Synthetic and Systems Biology Unit\\
Biological Research Centre\\
Hungarian Academy of Sciences\\
Szeged, Hungary\\

e-mail: jmolnar64@digikabel.hu, horvath.peter@brc.mta.hu}

\affil[$2$]{Institute for Molecular Medicine Finland (FIMM)\\
University of Helsinki\\
Helsinki, Finland}

\maketitle
\begin{abstract}
As the continuation of the contour mean calculation - designed for
averaging the manual delineations of 3D layer stack images - in this
paper, the most important equations: a) the reparameterization equations
to determine the minimizing diffeomorphism and b) the proper centroid
calculation for the surface mean calculation are presented. The chosen
representation space: Rescaled Position by Square root Normal (RPSN)
is a real valued vector space, invariant under the action of the reparameterization
group and the imposed $\mathbb{L}^{2}$ metric (used to define the
distance function) has well defined meaning: the sum of the central
second moments of the coordinate functions. For comparision purpose,
the reparameterization equations for elastic surface matching, using
the Square Root Normal Function (SRNF) are also provided. The reparameterization
equations for these cases have formal similarity, albeit the targeted
applications differ: SRNF representation suitable for shape analysis
purpose whereas RPSN is more fit for the cases where all contextual
information - including the relative translation between the constituent
surfaces - are to be retained (but the sake of theoretical completeness,
the possibility of the consistent relative displacement removal in
the RPSN case is also addressed).
\end{abstract}

\section{Inroduction}

Object delineation is an important annotation step to create training
data set for the supervised machine learning methods designed for
object segmentation. In cases wherever the object boundaries are not
definite (blurred, ambiguous), delineations performed by experts often
do not agree. A plausible approach to create meaningful annotation
samples is to accept the mean of many recommendations excluding some
outliers. This approach requires well defined, meaningful metrics
on the space of contours (2D) or surfaces (3D). A specifically designed
contour representation: the Rescaled Position by Square Velocity (RPSV):
$\mathbf{q}\left(t\right)=\mathbf{r}\left(t\right)\sqrt{\dot{\mathbf{r}}\left(t\right)}$
and a complete mathematical framework for the contour averaging (and
interpolation) problem were proposed and examined in the paper \cite{molnar2019elastic}.
RPSV can be considered as the mixing of two: the Kendall's 'landmark
points' \cite{Kendall84shapemanifolds} and the Square Root Velocity
Function (SRVF) \cite{joshi2007novel}\cite{10.1007/978-3-540-74198-5_30}\cite{srivastava2011shape}
representations. The description of the contours by preselected position
vector set (the landmark points) is the approach of the early shape
analysis techniques. The drawback of this representation is that the
results (\textit{e.g.} mean object) is dependent on the predetermined
sampling strategy of the landmark points (corresponding to fixed parameterization).
The continuous SRVF representation provides the necessary freedom
in the form of the optimal reparameterization of the contours. The
results are much much more intuitive. On the other hand, since the
description of the contours is velocity-vector based, the relative
displacement between the constituent contours cannot be retrieved
(insomuch as lost by derivation) - an important information for delineation
statistics. RPSV contour representation was designed to keep all contextual
information (including the relative displacenents of the constituent
contours) and covariant description. 

Recently, a new representation: the Square Root Normal Function (SRNF),
as the 3D generalization of the SRVF was introduced in \cite{jermyn2012elastic}
for elastic surface analysis. In this paper, the 3D generalization
- the 'Rescaled Position by Square root Normal' (RPSN) - of the RPSV
is introduced along with the reparameterization and the proper centroid
equations. Only these equations are presented, because the theoretical
results as the possibility of the pairwise determination of the optimal
parameterization system $\gamma_{i}\left(t\right)$, $\mathbf{q}_{i}\left(t\right)\rightarrow\mathbf{q}_{i}\left(\gamma_{i}\left(t\right)\right)$,
$i=0\ldots n-1$ wrt an arbitrarily chosen reference contour (\textit{e.g.}
$\gamma_{0}\left(t\right)\equiv t$) or the consistency of the mean
formula $\mathbf{q}\left(t\right)=\frac{1}{n}\underset{i}{\sum}\mathbf{q}_{i}\left(\gamma_{i}\left(t\right)\right)$
with the reparameterization equations, can be repeated and formally
transferred to 3D. The reparameterization equations are compared with
the corresponding equations of the SRNF. For the sake of theoretical
completeness, the possibility of the consistent relative displacement
removal (hence the generalization of the landmark point based approach
to the continuous case) is also shown. 

The structure of the paper is the following. Section \ref{sec:Elastic-surface-analysis}
summarizes the SRNF and RPSN frameworks. Section \ref{sec:Optimal-reparameterization-as}
is introduces to the reparameterization equations as variational problems,
section \ref{sec:Proper-centroid} is dedicated to the proper centroid
calculation (for RPSN) with reference to the possibility of relative
displacement removal between the constituent surfaces (\ref{sub:Removing-the-translation}).
Section \ref{sec:Conclusion} concludes the paper with discussion.
Appendices contain the important derivations.

\section{Elastic surface analysis frameworks\label{sec:Elastic-surface-analysis}}

We consider the primary space of the smooth surfaces embedded in $\mathbb{R}^{3}$
as the space of two-parameter coordinate function-triplets $X\left(u,v\right),\, Y\left(u,v\right),\, Z\left(u,v\right)$,
or - equivalently - two-parameter position vectors $\mathbf{S}\left(u,v\right)=X\mathbf{i}+Y\mathbf{j}+Z\mathbf{k}$
wrt some basis $\mathbf{i},\,\mathbf{j},\,\mathbf{k}$ (known as Gaussian
descriprion). The important derived quantities associated with this
description are: a) the local (covariant basis): $\mathbf{S}_{u}=\frac{\partial\mathbf{S}}{\partial u}$,
$\mathbf{S}_{v}=\frac{\partial\mathbf{S}}{\partial v}$, b) the normal
vector of the surface $\mathbf{N}=\mathbf{S}_{u}\times\mathbf{S}_{v}$,
c) its length $\left|\mathbf{N}\right|$ - which is also the square
root of the determinant of the metric tensor with components $g_{ik}=\mathbf{S}_{i}\times\mathbf{S}_{k},\, i,k\in\left[u,v\right]$,
$g=\det\left[g_{ik}\right]$,\textit{ i.e.} $\sqrt{g}=\left|\mathbf{N}\right|$
and d) the unit (paremeterization independent) normal vector $\mathbf{n}=\frac{\mathbf{N}}{\left|\mathbf{N}\right|}$.

Both representations SRNF: $\mathbf{P}\doteq\mathbf{n}\sqrt{\left|\mathbf{S}_{u}\times\mathbf{S}_{v}\right|}$
and RPSN: $\mathbf{Q}\doteq\mathbf{S}\sqrt{\left|\mathbf{S}_{u}\times\mathbf{S}_{v}\right|}$
can be considered as the ``change of coordinates'' in the space
of surfaces, albeit these mappings do not necessarily lead bijections
between the original and transformed coordinates. In the SRNF case,
any constant translation $\mathbf{d}$ between two surfaces $\mathbf{S}$
and $\tilde{\mathbf{S}}=\mathbf{S}+\mathbf{d}$ is removed by derivation,
so that in SRNF 'coordinates' the surfaces are determined only up
to an arbitrary translation, naturally representing the quotient space
of the surfaces wrt translation - a useful property for the shape
analysis.

These representetions have the distinguishing property that is the
distance functions defined between two points in the space of contours,
using the respective $\mathbb{L}^{2}$ norms:
\begin{equation}
d_{\mathbf{P}}^{2}=\iint\left(\mathbf{P}_{1}-\mathbf{P}_{2}\right)^{2}dudv\label{eq:srnf-dis2}
\end{equation}
and
\begin{equation}
d_{\mathbf{Q}}^{2}=\iint\left(\mathbf{Q}_{1}-\mathbf{Q}_{2}\right)^{2}dudv\label{eq:rpsn-dis2}
\end{equation}
are invariant under the action of the reparameterization group $\varGamma$
(in fact invariant quantities under the product group of the rotation
and reparameterization as well), because both represent parameterization
independent quantities: the surface area (SRNF) and the second central
moment (RPSN) of the surfaces. Among the consequences, a unique distance
over the equivalence class of the reparameterization group $\varGamma$
can be defined. Also, for multiple surface problems (\textit{e.g.}
averaging) any of the constituent surfaces can be chosen as reference
surface (wrt which the best reparameterizations of the other constituents
are to be determined). The transformation of the representations under
the effect of an arbitrary element $u\rightarrow\mu\left(u,v\right)$,
$\nu=\nu\left(u,v\right)$, $\left(\mu,\nu\right)\in\varGamma$ of
the reparameterization group are $\mathbf{P}\rightarrow\mathbf{P}\sqrt{\left|J_{\left(\mu,\nu\right)}\right|}$
and $\mathbf{Q}\rightarrow\mathbf{Q}\sqrt{\left|J_{\left(\mu,\nu\right)}\right|}$,
where the $\left|J_{\left(\mu,\nu\right)}\right|$ is the determinant
of the Jacobian $\left[J_{\left(\mu,\nu\right)}\right]=\left[\begin{array}{cc}
\frac{\partial\mu}{\partial u} & \frac{\partial\mu}{\partial v}\\
\frac{\partial\nu}{\partial u} & \frac{\partial\nu}{\partial v}
\end{array}\right]$ of the reparameterization. Reparameterization acting in the representation
space, require the updation of both the lengths and the directions
of the points of $\mathbf{P}$ or $\mathbf{Q}$. Alternatively however,
the distance minimization can be formulated directly in the surface
space (where only the direction of the points need to be updated),
the approach pursued in this paper.

\section{Optimal reparameterization as variational problem\label{sec:Optimal-reparameterization-as}}

Let $\mathbf{R}\left(u,v\right)$ with normal vector $\mathbf{M}=\mathbf{R}_{u}\times\mathbf{R}_{v}$
and unit normal $\mathbf{m}=\frac{\mathbf{M}}{\left|\mathbf{M}\right|}$
be the 'reference' surface and $\mathbf{S}\left(\mu,\nu\right)$ with
normal vectors $\mathbf{N}_{\left(\mu,\nu\right)}=\mathbf{S}_{\mu}\times\mathbf{S}_{\nu}$,
$\mathbf{n}=\frac{\mathbf{N}}{\left|\mathbf{N}\right|}$ another surface.
Let the $\mu\left(u,v\right)$, $\nu\left(u,v\right)$ function duplet
is the element of the reparameterization group $\varGamma$. We need
the surface $\mathbf{S}\left(u,v\right)=\mathbf{S}\left(\mu\left(u,v\right),\nu\left(u,v\right)\right)$
(with normal $\mathbf{N}=\mathbf{S}_{u}\times\mathbf{S}_{v}=\left|J_{\left(\mu,\nu\right)}\right|\mathbf{N}_{\left(\mu,\nu\right)}$)
to be optimally parameterized to the reference surface wrt the distance
function (\ref{eq:srnf-dis2}). Then the minimum distances between
surfaces can be formulated as variational problems using the 'direct
surface coordinates' in the space of surfaces instead of their representations.
In the case (\ref{eq:srnf-dis2}) the functional to be minimized is
\begin{eqnarray}
E\left(\mu,\nu\right) & = & \iint\left(\mathbf{m}\left(u,v\right)\sqrt{\left|\mathbf{R}_{u}\times\mathbf{R}_{v}\right|}-\mathbf{n}\left(u,v\right)\sqrt{\left|\mathbf{S}_{u}\times\mathbf{S}_{v}\right|}\right)^{2}dudv\nonumber \\
 & = & \iint\left(\mathbf{m}\left(u,v\right)\sqrt{\left|\mathbf{R}_{u}\times\mathbf{R}_{v}\right|}-\mathbf{n}\left(\mu,\nu\right)\sqrt{\left|J_{\left(\mu,\nu\right)}\right|\left|\mathbf{S}_{\mu}\times\mathbf{S}_{\nu}\right|}\right)^{2}dudv\label{eq:srnf-energy}
\end{eqnarray}
where $\left|J_{\left(\mu,\nu\right)}\right|$ stands for the determinant
of the Jacobian of the reparameterization:
\begin{eqnarray}
\left[J_{\left(\mu,\nu\right)}\right] & = & \left[\begin{array}{cc}
\frac{\partial\mu}{\partial u} & \frac{\partial\mu}{\partial v}\\
\frac{\partial\nu}{\partial u} & \frac{\partial\nu}{\partial v}
\end{array}\right],\nonumber \\
\left|J_{\left(\mu,\nu\right)}\right| & = & \frac{\partial\mu}{\partial u}\frac{\partial\nu}{\partial v}-\frac{\partial\mu}{\partial v}\frac{\partial\nu}{\partial u},
\end{eqnarray}
\textit{i.e.} the transformation between the normal vectors $\mathbf{N}=\mathbf{S}_{u}\times\mathbf{S}_{v}=\left|J_{\left(\mu,\nu\right)}\right|\mathbf{N}_{\left(\mu,\nu\right)}$
is used. The Euler-Lagrange equations associated with the minimization
problem (\ref{eq:srnf-energy}) (see \nameref{sec:Appendix-A}) are
\[
\left[\begin{array}{c}
\nabla_{\mu}E\\
\nabla_{\nu}E
\end{array}\right]=\sqrt{\left|\mathbf{M}\right|\left|\mathbf{N}\right|}\left[J_{\left(\mu,\nu\right)}\right]^{-1}\left[\begin{array}{c}
\mathbf{m}_{u}\cdot\mathbf{n}-\mathbf{m}\cdot\mathbf{n}_{u}+\frac{1}{2}\mathbf{m}\cdot\mathbf{n}\left(\Gamma_{u}^{\mathbf{R}}-\Gamma_{u}^{\mathbf{S}}\right)\\
\mathbf{m}_{v}\cdot\mathbf{n}-\mathbf{m}\cdot\mathbf{n}_{v}+\frac{1}{2}\mathbf{m}\cdot\mathbf{n}\left(\Gamma_{v}^{\mathbf{R}}-\Gamma_{v}^{\mathbf{S}}\right)
\end{array}\right],
\]
where $\Gamma_{w}^{\mathbf{R}}=\frac{\partial\ln\left|\mathbf{M}\right|}{\partial w}$,
$\Gamma_{w}^{\mathbf{S}}=\frac{\partial\ln\left|\mathbf{N}\right|}{\partial w}$,
$w\in\left[u,v\right]$ are the Christoffel divergences of the surfaces
$\mathbf{R}$ and $\mathbf{S}$ respectively. At around the identity
diffeomorphism these equations are simplified to the
\begin{eqnarray*}
\mathbf{m}_{u}\cdot\mathbf{n}-\mathbf{m}\cdot\mathbf{n}_{u}+\frac{1}{2}\mathbf{m}\cdot\mathbf{n}\left(\Gamma_{u}^{\mathbf{R}}-\Gamma_{u}^{\mathbf{S}}\right) & = & 0\\
\mathbf{m}_{v}\cdot\mathbf{n}-\mathbf{m}\cdot\mathbf{n}_{v}+\frac{1}{2}\mathbf{m}\cdot\mathbf{n}\left(\Gamma_{v}^{\mathbf{R}}-\Gamma_{v}^{\mathbf{S}}\right) & = & 0
\end{eqnarray*}
 ones that need to be solved for the new point positions on the fix-shaped
surface $\mathbf{S}$.

The associated functional to the problem (\ref{eq:rpsn-dis2}) is
\begin{eqnarray}
F\left(\mu,\nu\right) & = & \iint\left(\mathbf{R}\left(u,v\right)\sqrt{\left|\mathbf{R}_{u}\times\mathbf{R}_{v}\right|}-\mathbf{S}\left(u,v\right)\sqrt{\left|\mathbf{S}_{u}\times\mathbf{S}_{v}\right|}\right)^{2}dudv\nonumber \\
 & = & \iint\left(\mathbf{R}\left(u,v\right)\sqrt{\left|\mathbf{R}_{u}\times\mathbf{R}_{v}\right|}-\mathbf{S}\left(\mu,\nu\right)\sqrt{\left|J_{\left(\mu,\nu\right)}\right|\left|\mathbf{S}_{\mu}\times\mathbf{S}_{\nu}\right|}\right)^{2}dudv\label{eq:rpsn-energy}
\end{eqnarray}
and the associated Euler-Lagrange equations (see \nameref{sec:Appendix-B})
are:
\[
\left[\begin{array}{c}
\nabla_{\mu}F\\
\nabla_{\nu}F
\end{array}\right]=\sqrt{\left|\mathbf{M}\right|\left|\mathbf{N}\right|}\left[J_{\left(\mu,\nu\right)}\right]^{-1}\left[\begin{array}{c}
\mathbf{R}_{u}\cdot\mathbf{S}-\mathbf{R}\cdot\mathbf{S}_{u}+\frac{1}{2}\mathbf{R}\cdot\mathbf{S}\left(\Gamma_{u}^{\mathbf{R}}-\Gamma_{u}^{\mathbf{S}}\right)\\
\mathbf{R}_{v}\cdot\mathbf{S}-\mathbf{R}\cdot\mathbf{S}_{v}+\frac{1}{2}\mathbf{R}\cdot\mathbf{S}\left(\Gamma_{v}^{\mathbf{R}}-\Gamma_{v}^{\mathbf{S}}\right)
\end{array}\right]\,.
\]
At around the identity diffeomorphism these equations are simplified
to
\begin{eqnarray}
\mathbf{R}_{u}\cdot\mathbf{S}-\mathbf{R}\cdot\mathbf{S}_{u}+\frac{1}{2}\mathbf{R}\cdot\mathbf{S}\left(\Gamma_{u}^{\mathbf{R}}-\Gamma_{u}^{\mathbf{S}}\right) & = & 0\nonumber \\
\mathbf{R}_{v}\cdot\mathbf{S}-\mathbf{R}\cdot\mathbf{S}_{v}+\frac{1}{2}\mathbf{R}\cdot\mathbf{S}\left(\Gamma_{v}^{\mathbf{R}}-\Gamma_{v}^{\mathbf{S}}\right) & = & 0\,.\label{eq:RPSN-E-L_at_identity}
\end{eqnarray}
As in the two dimensional (contour) case, it can be shown that for
the mean surface problem:
\begin{equation}
\underset{\mu_{i},\nu_{i}}{\min}\stackrel[i=0]{N-1}{\sum}\iint\left(\mathbf{Q}-\mathbf{Q_{i}}\right)^{2}dudv,\label{eq:rpsn-system-energy}
\end{equation}
the solution is $\mathbf{Q}\left(u,v\right)=\frac{1}{N}\stackrel[i=0]{N-1}{\sum}\mathbf{Q_{i}}\left(u,v\right)$,
where the constituents $\left(\mathbf{Q_{i}}\right)$ are all optimally
parameterized wrt a freely chosen reference surface (say $\mathbf{Q_{0}}$),
\textit{i.e.} $\mathbf{Q_{i}}\left(u,v\right)=\mathbf{Q_{i}}\left(\mu_{i}\left(u,v\right),\nu_{i}\left(u,v\right)\right)$,
$i=1...N-1$, where $\mu_{i}\left(u,v\right),\:\nu_{i}\left(u,v\right)$
are the optimal diffeomorphisms - the solution of the Euler-Lagrange
equation systems (\ref{eq:RPSN-E-L_at_identity}) - determined pairwise
between $\mathbf{Q_{0}}=\mathbf{R}\left(u,v\right)\sqrt{\left|\mathbf{R}_{u}\times\mathbf{R}_{v}\right|}$
and 
\begin{gather*}
\mathbf{Q_{i}}=\mathbf{S_{i}}\left(\mu_{i},\nu_{i}\right)\sqrt{\left|J_{\left(\mu_{i},\nu_{i}\right)}\right|\left|\mathbf{S_{i}}_{\mu_{i}}\left(\mu_{i},\nu_{i}\right)\times\mathbf{S_{i}}_{v_{i}}\left(\mu_{i},\nu_{i}\right)\right|}\\
\qquad\qquad\qquad\qquad\qquad=\mathbf{S_{i}}\left(u,v\right)\sqrt{\left|\mathbf{S_{i}}_{u}\left(u,v\right)\times\mathbf{S_{i}}_{v}\left(u,v\right)\right|}
\end{gather*}
since $\mathbf{S_{i}}\left(u,v\right)=\mathbf{S_{i}}\left(\mu_{i}\left(u,v\right),\nu_{i}\left(u,v\right)\right)$.

\section{Proper centroid\label{sec:Proper-centroid}}

In the RPSN case, the Euler-Lagrange equations for the optimal reparameterization
system can be solved for surfaces (that are given as position vectors
$\mathbf{S}=X\mathbf{i}+Y\mathbf{j}+Z\mathbf{k}$) wrt any point in
the space designated as origin, but the result is dependent on the
choice of the origin. However there is an optimal choice of the origin
compatible with the minimization problem (\ref{eq:rpsn-system-energy}).
The new energy comprising the position of this optimal cetroid - wrt
the current origin - can be formulated as double-minimization problem:
\begin{equation}
\underset{\mu_{i},\nu_{i};\delta\mathbf{D}}{\min}\stackrel[i=0]{N-1}{\sum}\iint\left[\left(\mathbf{S}-\delta\mathbf{D}\right)\sqrt{\left|\mathbf{N}\right|}-\left(\mathbf{S_{i}}-\delta\mathbf{D}\right)\sqrt{\left|\mathbf{N_{i}}\right|}\right]^{2}dudv,\label{eq:rpsn-system-energy-1}
\end{equation}
where the new representations (wrt the displaced origin) are $\mathbf{\mathbf{Q}^{\delta D}=\left(\mathbf{S}-\delta\mathbf{D}\right)}\sqrt{\left|\mathbf{N}\right|}=\frac{1}{N}\stackrel[i=0]{N-1}{\sum}\mathbf{Q_{i}^{\delta D}}$,
$\mathbf{\frac{1}{N}\stackrel[i=0]{N-1}{\sum}\mathbf{Q_{i}^{\delta D}}=\left(\mathbf{S_{i}}-\delta\mathbf{D}\right)}\frac{1}{N}\stackrel[i=0]{N-1}{\sum}\mathbf{S_{i}}\sqrt{\left|\mathbf{N_{i}}\right|}$
(the reference contour $\mathbf{\mathbf{R}=S_{0}}$, $\mathbf{N}=\mathbf{S}_{u}\times\mathbf{S}_{v}$,
$\mathbf{N_{i}}=\mathbf{S_{i}}_{u}\times\mathbf{S_{i}}_{v}$). Notation
$\delta\mathbf{D}$ indicates that the new origin is to be calculated
wrt fixed set of optimally (precedingly) determined reparameterization.
The minimization wrt the position of the new origin leads to simple
extreme value problem: 
\begin{equation}
\frac{\partial}{\partial\delta\mathbf{D}}\stackrel[i=0]{N-1}{\sum}\iint\left[\left(\mathbf{S}-\delta\mathbf{D}\right)\sqrt{\left|\mathbf{N}\right|}-\left(\mathbf{S_{i}}-\delta\mathbf{D}\right)\sqrt{\left|\mathbf{N_{i}}\right|}\right]^{2}dudv\doteq0,\label{eq:rpsn-system-derivative}
\end{equation}
the solution is 
\begin{equation}
\delta\mathbf{D}=\frac{\stackrel[i=0]{N-1}{\sum}\iint\left(\mathbf{S_{i}}\sqrt{\left|\mathbf{N_{i}}\right|}-\mathbf{S}\sqrt{\left|\mathbf{N}\right|}\right)\sqrt{\left|\mathbf{N_{i}}\right|}dudv}{\stackrel[i=0]{N-1}{\sum}\iint\left(\sqrt{\left|\mathbf{N_{i}}\right|}-\sqrt{\left|\mathbf{N}\right|}\right)^{2}dudv}=0\,.\label{eq:proper-centroid}
\end{equation}
The optimal set of diffeomorphism and the position of the proper centroid
are interrelated: the new centroid position involves a new set of
optimal reparameterization functions $\mu_{i}\left(u,v\right),\nu_{i}\left(u,v\right)$,
$i=1...N-1$, which in turn moves the optimal centroid further away
from its ad hoc initial position. Therefore the double optimization
problem (\ref{eq:rpsn-system-energy-1}) can be solved iteratively.

\subsection{Removing translations\label{sub:Removing-the-translation}}

The RPSN representation is intentionally designed to retain the relative
displacement information between the constituent surfaces, nevertheless
it could be used for shape analysis purpose by appropriately removing
the relative displacements between the constituents. Assume the double
optimization problem (\ref{eq:rpsn-system-energy-1}) is solved .
Then the minimization problem: 
\begin{equation}
\underset{\delta\mathbf{D_{i}}}{\min}\stackrel[i=0]{N-1}{\sum}\iint\left[\mathbf{S}\sqrt{\left|\mathbf{N}\right|}-\left(\mathbf{S_{i}}-\delta\mathbf{D_{i}}\right)\sqrt{\left|\mathbf{N_{i}}\right|}\right]^{2}dudv,\label{eq:rpsn-system-energy-1-1}
\end{equation}
can be used to remove the displacement updating the original positions
of the constituents $\mathbf{S_{i}}$ to$\mathbf{S_{i}}-\delta\mathbf{D_{i}}$.
From the extreme value problems: 
\begin{equation}
\frac{\partial}{\partial\delta\mathbf{D}}\stackrel[i=0]{N-1}{\sum}\iint\left[\mathbf{S}\sqrt{\left|\mathbf{N}\right|}-\left(\mathbf{S_{i}}-\delta\mathbf{D_{i}}\right)\sqrt{\left|\mathbf{N_{i}}\right|}\right]^{2}dudv\doteq0,\:\mathbf{i}=0...N-1\label{eq:rpsn-system-derivative-1}
\end{equation}
and the solutions are
\begin{equation}
\delta\mathbf{D_{i}}=\frac{\iint\left(\mathbf{S_{i}}\sqrt{\left|\mathbf{N_{i}}\right|}-\mathbf{S}\sqrt{\left|\mathbf{N}\right|}\right)\sqrt{\left|\mathbf{N_{i}}\right|}dudv}{\iint\left|\mathbf{N_{i}}\right|dudv}=0.\label{eq:removing-displacement}
\end{equation}
(note the denominator is the surface area of the $i$-th constituent
surface. 

As the optimal set of diffeomorphism and the position of the proper
centroid are interrelated, so the displacements calculated by (\ref{eq:removing-displacement}).
Iterative solution is possible. Since the quantities that need to
be calculated for the proper centroid (\ref{eq:proper-centroid})
and the displacement removal equations (\ref{eq:removing-displacement})
are largely overlapped, $\delta\mathbf{D}$ and all $\delta\mathbf{D_{i}}$,
$\mathbf{i}=0...N-1$ should expediently be calculated in the same
step. In this case, however the problem of the convergence remains
an open question.

\section{Discussion\label{sec:Conclusion}}

In this paper an elastic surface mean determination method was presented.
The mean surface is calculated from a set of surfaces in a way that
all visible information (relative placement, rotation, scale) are
retained. At the same time - borrowing the idea from the state of
the art shape analysis methods - the parameterization of the surfaces
is relaxed. The chosen contour representation (RPSN) and the imposed
$\mathbb{L}^{2}$ metric forms a Hilbert space of the contour representations.
The metric is chosen to be invariant wrt the reparameterization, the
distance function based on it has well defined meaning, the (sum of)
the second moment of the surfaces. The mean surface calculation is
performed in the quotient space space of surfaces modulo reparameterization
group and could be formulated as a double optimization problem: a
variational for the system of the optimal parameterization and an
extreme value problem for the proper centroid identification. 

The work is the direct extension of the 2D contour mean calculation
using the representation RPSV (Rescaled Position by Square Velocity)
inroduced in the paper \cite{molnar2019elastic}.

\pagebreak{}

\part*{Appendices}

In the appendices, the optimal reparameterization for the SRNF (\nameref{sec:Appendix-A})
and RPSN (\nameref{sec:Appendix-B}) are derived. The common notations
and basic formulae used in these appendices are introduced below.

The notations wrt the original parameterization $\left(u,v\right)$
are as follows. The covariant basis vectors of the surface $\mathbf{S}=\mathbf{S}\left(u,v\right)$
are denoted by $\mathbf{S}_{u}=\frac{\partial\mathbf{S}}{\partial u}$,
$\mathbf{S}_{v}=\frac{\partial\mathbf{S}}{\partial v}$. The contravariant
basis vectors $\mathbf{S}^{u}$, $\mathbf{S}^{v}$ are defined such
that their dot (scalar) products with the coordinate basis vectors
are $\mathbf{S}^{i}\cdot\mathbf{S}_{k}=\delta_{k}^{i}$, $i,k\in\left\{ u,v\right\} $.
The normal vector of the surface is denoted by $\mathbf{N}=\mathbf{S}_{u}\times\mathbf{S}_{v}$.
$\left|\mathbf{N}\right|$ stands for its length, the unit (hence
paremeterization-independent) normal vector is denoted by $\mathbf{n}$
$\left(\mathbf{n}=\frac{\mathbf{N}}{\left|\mathbf{N}\right|}\right)$.
The direct (dyadic) product of two vectors $\mathbf{u}$, $\mathbf{v}$
is defined such that its scalar products (contractions) with the vector
$\mathbf{w}$ become $\left(\mathbf{uv}\right)\cdot\mathbf{w}=\left(\mathbf{v}\cdot\mathbf{w}\right)\mathbf{u}$
and $\mathbf{w}\cdot\left(\mathbf{uv}\right)=\left(\mathbf{u}\cdot\mathbf{w}\right)\mathbf{v}$.
Metric and inverse metric tensor components are defined as $g_{ik}=\mathbf{S}_{i}\cdot\mathbf{S}_{k}$
and $g^{ik}=\mathbf{S}^{i}\cdot\mathbf{S}^{k}$, $i,k\in\left\{ u,v\right\} $
respectively. The determinant of the metric tensor is denoted by $g$
$\left(g=\det\left[g_{ik}\right]\right)$, its square root (the 'metric')
$\sqrt{g}=\left|\mathbf{N}\right|$ is used to define surface area
element as $dS=\sqrt{g}dudv$. Christoffel symbols (connection components)
for embedded surfaces can be defined as 
\begin{gather}
\Gamma_{ik}^{l}=\mathbf{S}^{l}\cdot\mathbf{S}_{ik}\nonumber \\
i,\, k,\, l\in\left\{ u,v\right\} ,\label{eq:ChristoffelS}
\end{gather}
where vectors $\mathbf{S}_{ik}$ are the second partial derivative
of the position vector $\mathbf{S}$ . It can be seen by simple substitution
that
\begin{eqnarray}
\mathbf{S}^{u} & = & \frac{1}{\left|\mathbf{S}_{u}\times\mathbf{S}_{v}\right|}\mathbf{S}_{v}\times\mathbf{n}\nonumber \\
\mathbf{S}^{v} & = & \frac{1}{\left|\mathbf{S}_{u}\times\mathbf{S}_{v}\right|}\mathbf{n}\times\mathbf{S}_{u}\,.\label{eq:contra-by-co}
\end{eqnarray}
Using (\ref{eq:ChristoffelS}), (\ref{eq:contra-by-co}), the partial
derivatives of the logarithm of the metric $\left|\mathbf{N}\right|$
become the Christoffel 'divergences':
\begin{gather}
\frac{\partial\ln\left|\mathbf{N}\right|}{\partial u}=\frac{1}{\left|\mathbf{N}\right|}\mathbf{n}\cdot\left(\mathbf{S}_{uu}\times\mathbf{S}_{v}+\mathbf{S}_{u}\times\mathbf{S}_{uv}\right)\nonumber \\
=\mathbf{S}_{uu}\cdot\mathbf{S}^{u}+\mathbf{S}_{uv}\cdot\mathbf{S}^{v}\nonumber \\
=\Gamma_{uu}^{u}+\Gamma_{vu}^{v}\label{eq:Christdivu}\\
\frac{\partial\ln\left|\mathbf{N}\right|}{\partial v}=\frac{1}{\left|\mathbf{N}\right|}\mathbf{n}\cdot\left(\mathbf{S}_{uv}\times\mathbf{S}_{v}+\mathbf{S}_{u}\times\mathbf{S}_{vv}\right)\nonumber \\
=\mathbf{S}_{uv}\cdot\mathbf{S}^{u}+\mathbf{S}_{vv}\cdot\mathbf{S}^{v}\nonumber \\
=\Gamma_{uv}^{u}+\Gamma_{vv}^{v}\label{eq:Christdivv}
\end{gather}
From identity $\mathbf{n}\cdot\mathbf{S}_{k}\equiv0$ 
\begin{gather}
\mathbf{n}_{i}\cdot\mathbf{S}_{k}=-\mathbf{n}\cdot\mathbf{S}_{ik}\nonumber \\
i,k\in\left\{ u,v\right\} \,.\label{eq:to-secfundform}
\end{gather}
Any vector $\mathbf{w}$ can be decomposed in the local basis $\mathbf{S}_{u},\,\mathbf{S}_{v},\,\mathbf{n}$
as $\mathbf{w}=\left(\mathbf{w}\cdot\mathbf{S}^{u}\right)\mathbf{S}_{u}+\left(\mathbf{w}\cdot\mathbf{S}^{v}\right)\mathbf{S}_{v}+\left(\mathbf{w}\cdot\mathbf{n}\right)\mathbf{n}$
or alternativeli in the contravariant basis as $\mathbf{w}=\left(\mathbf{w}\cdot\mathbf{S}_{u}\right)\mathbf{S}^{u}+\left(\mathbf{w}\cdot\mathbf{S}_{v}\right)\mathbf{S}^{v}+\left(\mathbf{w}\cdot\mathbf{n}\right)\mathbf{n}$,
it follows that the decomposition of the identity tensor ($\mathbf{I}:\,\mathbf{I}\cdot\mathbf{w}=\mathbf{w}\cdot\mathbf{I}\equiv\mathbf{w}$)
can be given in two ways:
\begin{eqnarray*}
\mathbf{I} & = & \mathbf{S}^{u}\mathbf{S}_{u}+\mathbf{S}^{v}\mathbf{S}_{v}+\mathbf{nn}\\
\mathbf{I} & = & \mathbf{S}_{u}\mathbf{S}^{u}+\mathbf{S}_{v}\mathbf{S}^{v}+\mathbf{nn\,.}
\end{eqnarray*}
Note that the partial derivatives of the unit normal vector $\mathbf{n}_{u}$
and $\mathbf{n}_{v}$ are the elements of the tangent space hence
can be decomposed such that $\mathbf{n}_{i}=\left(\mathbf{n}_{i}\cdot\mathbf{S}^{u}\right)\mathbf{S}_{u}+\left(\mathbf{n}_{i}\cdot\mathbf{S}^{v}\right)\mathbf{S}_{v}$. 

The partial derivatives of the SRNF representation $\mathbf{P}=\mathbf{n}\sqrt{\left|\mathbf{S}_{u}\times\mathbf{S}_{v}\right|}$:
\begin{eqnarray}
\mathbf{P}_{u} & = & \sqrt{\left|\mathbf{N}\right|}\mathbf{n}_{u}+\mathbf{n}\frac{\mathbf{n}}{2\sqrt{\left|\mathbf{N}\right|}}\cdot\left(\mathbf{S}_{uu}\times\mathbf{S}_{v}+\mathbf{S}_{u}\times\mathbf{S}_{uv}\right)\nonumber \\
 & = & \sqrt{\left|\mathbf{N}\right|}\mathbf{n}_{u}+\mathbf{n}\frac{1}{2\sqrt{\left|\mathbf{N}\right|}}\left[\mathbf{S}_{uu}\cdot\left(\mathbf{S}_{v}\times\mathbf{n}\right)+\mathbf{S}_{uv}\cdot\left(\mathbf{n}\times\mathbf{S}_{u}\right)\right]\label{eq:SRNF_u}\\
 & = & \sqrt{\left|\mathbf{N}\right|}\left[\mathbf{n}_{u}+\frac{\mathbf{n}}{2}\left(\Gamma_{uu}^{u}+\Gamma_{vu}^{v}\right)\right],\nonumber 
\end{eqnarray}
similarly 
\begin{eqnarray}
\mathbf{P}_{v} & = & \sqrt{\left|\mathbf{N}\right|}\left[\mathbf{n}_{v}+\frac{\mathbf{n}}{2}\left(\Gamma_{vu}^{v}+\Gamma_{vv}^{v}\right)\right]\,.\label{eq:SRNF_v}
\end{eqnarray}

The partial derivatives of the RPSN representation $\mathbf{Q}=\mathbf{S}\sqrt{\left|\mathbf{S}_{u}\times\mathbf{S}_{v}\right|}$:
\begin{eqnarray}
\mathbf{Q}_{u} & = & \sqrt{\left|\mathbf{N}\right|}\mathbf{S}_{u}+\mathbf{S}\frac{\mathbf{n}}{2\sqrt{\left|\mathbf{N}\right|}}\cdot\left(\mathbf{S}_{uu}\times\mathbf{S}_{v}+\mathbf{S}_{u}\times\mathbf{S}_{uv}\right)\nonumber \\
 & = & \sqrt{\left|\mathbf{N}\right|}\mathbf{S}_{u}+\mathbf{S}\frac{1}{2\sqrt{\left|\mathbf{N}\right|}}\left[\mathbf{S}_{uu}\cdot\left(\mathbf{S}_{v}\times\mathbf{n}\right)+\mathbf{S}_{uv}\cdot\left(\mathbf{n}\times\mathbf{S}_{u}\right)\right]\label{eq:RPSN_u}\\
 & = & \sqrt{\left|\mathbf{N}\right|}\left[\mathbf{S}_{u}+\frac{\mathbf{S}}{2}\left(\Gamma_{uu}^{u}+\Gamma_{vu}^{v}\right)\right]\nonumber 
\end{eqnarray}
similarly 
\begin{eqnarray}
\mathbf{Q}_{v} & = & \sqrt{\left|\mathbf{N}\right|}\left[\mathbf{S}_{v}+\frac{\mathbf{S}}{2}\left(\Gamma_{vu}^{v}+\Gamma_{vv}^{v}\right)\right]\,.\label{eq:RPSN_v}
\end{eqnarray}

The transformation of the normal $\mathbf{N}=\mathbf{S}_{u}\times\mathbf{S}_{v}\rightarrow\mathbf{N}_{\left(\mu,\nu\right)}=\mathbf{S}_{\mu}\times\mathbf{S}_{\nu}$
can be expressed as 
\begin{gather}
\mathbf{S}_{u}\times\mathbf{S}_{v}=\left(\mathbf{S}_{\mu}\frac{\partial\mu}{\partial u}+\mathbf{S}_{\nu}\frac{\partial\nu}{\partial u}\right)\times\left(\mathbf{S}_{\mu}\frac{\partial\mu}{\partial v}+\mathbf{S}_{\nu}\frac{\partial\nu}{\partial v}\right)\nonumber \\
=\left(\frac{\partial\mu}{\partial u}\frac{\partial\nu}{\partial v}-\frac{\partial\mu}{\partial v}\frac{\partial\nu}{\partial u}\right)\mathbf{S}_{\mu}\times\mathbf{S}_{\nu}\label{eq:metric-transformation}\\
=\left|J_{\left(\mu,\nu\right)}\right|\mathbf{S}_{\mu}\times\mathbf{S}_{\nu}\nonumber 
\end{gather}
where 
\begin{eqnarray}
\left|J_{\left(\mu,\nu\right)}\right| & = & \frac{\partial\mu}{\partial u}\frac{\partial\nu}{\partial v}-\frac{\partial\mu}{\partial v}\frac{\partial\nu}{\partial u}\,.\label{eq:reparameterization-Jacobiandeterminant}
\end{eqnarray}
 is the determinant of the Jacobian 
\begin{eqnarray}
\left[J_{\left(\mu,\nu\right)}\right] & = & \left[\begin{array}{cc}
\frac{\partial\mu}{\partial u} & \frac{\partial\mu}{\partial v}\\
\frac{\partial\nu}{\partial u} & \frac{\partial\nu}{\partial v}
\end{array}\right]
\end{eqnarray}
of the reparameterization. The reparameterization $\left(\mu\left(u,v\right),\nu\left(u,v\right)\right)$
considered feasible iff its determinant (\ref{eq:reparameterization-Jacobiandeterminant})
is not negative for any values of $\left(u,v\right)$ and zero only
in isolated points. The inverse of the Jacobian is:
\begin{eqnarray}
\left[J_{\left(\mu,\nu\right)}\right]^{-1} & =\left[\begin{array}{cc}
\frac{\partial u}{\partial\mu} & \frac{\partial u}{\partial\nu}\\
\frac{\partial v}{\partial\mu} & \frac{\partial v}{\partial\nu}
\end{array}\right]=\frac{1}{\left|J_{\left(\mu,\nu\right)}\right|} & \left[\begin{array}{cc}
\frac{\partial\nu}{\partial v} & -\frac{\partial\mu}{\partial v}\\
-\frac{\partial\nu}{\partial u} & \frac{\partial\mu}{\partial u}
\end{array}\right]\,.\label{eq:inverse-Jacobian}
\end{eqnarray}
Using (\nameref{eq:metric-transformation}) and (\nameref{eq:inverse-Jacobian}),
the transformation of the Christoffel divergences become:
\begin{gather*}
\frac{\partial\ln\left|\mathbf{S}_{\mu}\times\mathbf{S}_{\nu}\right|}{\partial\mu}=\qquad\qquad\qquad\qquad\qquad\qquad\qquad\qquad\\
\frac{\partial\ln\frac{\left|\mathbf{S}_{u}\times\mathbf{S}_{v}\right|}{\left|J_{\left(\mu,\nu\right)}\right|}}{\partial u}\frac{\partial u}{\partial\mu}+\frac{\partial\ln\frac{\left|\mathbf{S}_{u}\times\mathbf{S}_{v}\right|}{\left|J_{\left(\mu,\nu\right)}\right|}}{\partial v}\frac{\partial v}{\partial\mu}=\\
\frac{\frac{\partial\nu}{\partial v}}{\left|J_{\left(\mu,\nu\right)}\right|}\left[-\frac{\left|J_{\left(\mu,\nu\right)}\right|_{u}}{\left|J_{\left(\mu,\nu\right)}\right|}+\Gamma_{uu}^{u}+\Gamma_{vu}^{v}\right]-\frac{\frac{\partial\nu}{\partial u}}{\left|J_{\left(\mu,\nu\right)}\right|}\left[-\frac{\left|J_{\left(\mu,\nu\right)}\right|_{v}}{\left|J_{\left(\mu,\nu\right)}\right|}+\Gamma_{uv}^{u}+\Gamma_{vv}^{v}\right]\\
\Rightarrow
\end{gather*}
\begin{gather}
\Gamma_{\mu\mu}^{\mu}+\Gamma_{v\mu}^{v}=\nonumber \\
\frac{\frac{\partial\nu}{\partial v}}{\left|J_{\left(\mu,\nu\right)}\right|}\left[-\frac{\left|J_{\left(\mu,\nu\right)}\right|_{u}}{\left|J_{\left(\mu,\nu\right)}\right|}+\Gamma_{uu}^{u}+\Gamma_{vu}^{v}\right]-\frac{\frac{\partial\nu}{\partial u}}{\left|J_{\left(\mu,\nu\right)}\right|}\left[-\frac{\left|J_{\left(\mu,\nu\right)}\right|_{v}}{\left|J_{\left(\mu,\nu\right)}\right|}+\Gamma_{uv}^{u}+\Gamma_{vv}^{v}\right],\label{eq:Christoffel-divergence-transform-mu}\\
\Gamma_{\mu\nu}^{\mu}+\Gamma_{v\nu}^{v}=...=\nonumber \\
-\frac{\frac{\partial\mu}{\partial v}}{\left|J_{\left(\mu,\nu\right)}\right|}\left[-\frac{\left|J_{\left(\mu,\nu\right)}\right|_{u}}{\left|J_{\left(\mu,\nu\right)}\right|}+\Gamma_{uu}^{u}+\Gamma_{vu}^{v}\right]+\frac{\frac{\partial\mu}{\partial u}}{\left|J_{\left(\mu,\nu\right)}\right|}\left[-\frac{\left|J_{\left(\mu,\nu\right)}\right|_{v}}{\left|J_{\left(\mu,\nu\right)}\right|}+\Gamma_{uv}^{u}+\Gamma_{vv}^{v}\right].\label{eq:Christoffel-divergence-transform-nu}
\end{gather}
Note that one can apply the transformation rule of the Christoffel
symbols directly to deduce results (\nameref{eq:Christoffel-divergence-transform-mu}),
(\nameref{eq:Christoffel-divergence-transform-nu}).

\section*{Appendix A\label{sec:Appendix-A}}

The optimal parameterization for the SRNF problem between surfaces
$\mathbf{R}\left(u,v\right)$ (reference) and $\mathbf{S}\left(u,v\right)=\mathbf{S}\left(\mu\left(u,v\right),\nu\left(u,v\right)\right)$
can be formulated as

\begin{equation}
\underset{\mu,\nu}{\min}\iint\left(\mathbf{m}\left(u,v\right)\sqrt{\left|\mathbf{R}_{u}\times\mathbf{R}_{v}\right|}-\mathbf{n}\left(\mu,\nu\right)\sqrt{\left|J_{\left(\mu,\nu\right)}\right|\left|\mathbf{S}_{\mu}\times\mathbf{S}_{\nu}\right|}\right)^{2}dudv\label{eq:srnf-min1}
\end{equation}
$\left(\left|\mathbf{S}_{u}\times\mathbf{S}_{v}\right|=\left|J_{\left(\mu,\nu\right)}\right|\left|\mathbf{S}_{\mu}\times\mathbf{S}_{\nu}\right|,\:\left|J_{\left(\mu,\nu\right)}\right|=\frac{\partial\mu}{\partial u}\frac{\partial\nu}{\partial v}-\frac{\partial\mu}{\partial v}\frac{\partial\nu}{\partial u}\right)$.
Since the square of the terms in the parentheses $\iint\left|\mathbf{R}_{u}\times\mathbf{R}_{v}\right|dudv$,
$\iint\left|\mathbf{S}_{u}\times\mathbf{S}_{v}\right|dudv$ express
the surface area of the surfaces $\mathbf{R}$ and $\mathbf{S}$,
problem (\ref{eq:srnf-min1}) is equivalent with
\begin{gather}
\underset{\mu,\nu}{\min}-\iint\mathbf{m}\left(u,v\right)\cdot\mathbf{n}\left(\mu,\nu\right)\sqrt{\left|\mathbf{R}_{u}\times\mathbf{R}_{v}\right|}\sqrt{\left|J_{\left(\mu,\nu\right)}\right|\left|\mathbf{S}_{\mu}\times\mathbf{S}_{\nu}\right|}dudv\,.\label{eq:srnf-min2}
\end{gather}
Its Lagrangian is
\begin{gather}
L=-\sqrt{\left|J_{\left(\mu,\nu\right)}\right|}\sqrt{\left|\mathbf{M}\right|}\mathbf{m}\cdot\mathbf{n}\sqrt{\left|\mathbf{S}_{\mu}\times\mathbf{S}_{\nu}\right|}\label{eq:srnf-Lagrangian}
\end{gather}
$\left(\mathbf{M}=\mathbf{R}_{u}\times\mathbf{R}_{v},\:\mathbf{m}=\frac{\mathbf{M}}{\left|\mathbf{M}\right|}\right)$.
From the latter 
\begin{gather}
\frac{\partial L}{\partial\mu}=-\sqrt{\left|J_{\left(\mu,\nu\right)}\right|}\sqrt{\left|\mathbf{M}\right|}\mathbf{m}\cdot\left[\sqrt{\left|\mathbf{S}_{\mu}\times\mathbf{S}_{\nu}\right|}\mathbf{n}_{\mu}+\mathbf{n}\frac{\mathbf{n}\cdot\left(\mathbf{S}_{\mu\mu}\times\mathbf{S}_{\nu}+\mathbf{S}_{\mu}\times\mathbf{S}_{\mu\nu}\right)}{2\sqrt{\left|\mathbf{S}_{\mu}\times\mathbf{S}_{\nu}\right|}}\right]\nonumber \\
=-\sqrt{\left|\mathbf{M}\right|\left|\mathbf{N}\right|}\mathbf{m}\cdot\left[\mathbf{n}_{\mu}+\mathbf{n}\frac{\mathbf{S}_{\mu\mu}\cdot\left(\mathbf{S}_{\nu}\times\mathbf{n}\right)+\mathbf{S}_{\mu\nu}\cdot\left(\mathbf{n}\times\mathbf{S}_{\mu}\right)}{2\left|\mathbf{S}_{\mu}\times\mathbf{S}_{\nu}\right|}\right]\label{eq:dL/dmu0}
\end{gather}
where notation $\left|\mathbf{N}\right|$ is exclusively used to denote
the normal vector of the surface wrt parameterization $\left(u,v\right)$,
\textit{i.e.} $\mathbf{N}=\mathbf{S}_{u}\times\mathbf{S}_{v}=\left|J_{\left(\mu,\nu\right)}\right|\mathbf{S}_{\mu}\times\mathbf{S}_{\nu}$
(see Eq. (\nameref{eq:metric-transformation})). Using the identity
(\ref{eq:contra-by-co}) and the definitions of the Christoffel symbols
(\ref{eq:ChristoffelS}), (\ref{eq:dL/dmu0}) can be rearranged as:
\begin{gather}
\frac{\partial L}{\partial\mu}=-\sqrt{\left|\mathbf{M}\right|\left|\mathbf{N}\right|}\left(\mathbf{m}\cdot\mathbf{n}_{\mu}+\frac{1}{2}\mathbf{m}\cdot\mathbf{n}\Gamma_{\mu}^{\mathbf{S}}\right)\label{eq:dL/dmu}
\end{gather}
and similarly
\begin{gather}
\frac{\partial L}{\partial\nu}=-\sqrt{\left|\mathbf{M}\right|\left|\mathbf{N}\right|}\left(\mathbf{m}\cdot\mathbf{n}_{\nu}+\frac{1}{2}\mathbf{m}\cdot\mathbf{n}\Gamma_{\nu}^{\mathbf{S}}\right),\label{eq:dL/dnu}
\end{gather}
where the
\begin{eqnarray}
\Gamma_{\mu}^{\mathbf{S}} & = & \mathbf{S}_{\mu\mu}\cdot\mathbf{S}^{\mu}+\mathbf{S}_{\mu\nu}\cdot\mathbf{S}^{\nu}=\frac{\partial\ln\left|\mathbf{N}\right|}{\partial\mu}\nonumber \\
\Gamma_{\nu}^{\mathbf{S}} & = & \mathbf{S}_{\mu\nu}\cdot\mathbf{S}^{\mu}+\mathbf{S}_{\nu\nu}\cdot\mathbf{S}^{\nu}=\frac{\partial\ln\left|\mathbf{N}\right|}{\partial\nu}\label{eq:Christoffel-divergences}
\end{eqnarray}
are the Christoffel divergences , \textit{i.e.} the relative change
of the metric $\left|\mathbf{N}\right|$ wrt the variables $\mu,\,\nu$.
Using (\ref{eq:inverse-Jacobian}) (\ref{eq:Christoffel-divergence-transform-mu}),
equations (\ref{eq:dL/dmu}) expressed with the original coordinates
$\left(u,v\right)$ can be rearranged as
\begin{gather}
\frac{2\left|J_{\left(\mu,\nu\right)}\right|}{\sqrt{\left|\mathbf{M}\right|\left|\mathbf{N}\right|}}\frac{\partial L}{\partial\mu}=-2\mathbf{m}\cdot\left(\frac{\partial\nu}{\partial v}\mathbf{n}_{u}-\frac{\partial\nu}{\partial u}\mathbf{n}_{v}\right)\qquad\qquad\qquad\qquad\qquad\qquad\qquad\nonumber \\
\qquad\qquad\qquad\qquad\qquad-\mathbf{m}\cdot\mathbf{n}\frac{\partial\nu}{\partial v}\left[-\frac{\left|J_{\left(\mu,\nu\right)}\right|_{u}}{\left|J_{\left(\mu,\nu\right)}\right|}+\Gamma_{u}^{\mathbf{S}}\right]+\mathbf{m}\cdot\mathbf{n}\frac{\partial\nu}{\partial u}\left[-\frac{\left|J_{\left(\mu,\nu\right)}\right|_{v}}{\left|J_{\left(\mu,\nu\right)}\right|}+\Gamma_{v}^{\mathbf{S}}\right]=\nonumber \\
\frac{\partial\nu}{\partial v}\left[-2\mathbf{m}\cdot\mathbf{n}_{u}-\mathbf{m}\cdot\mathbf{n}\left(-\frac{\left|J_{\left(\mu,\nu\right)}\right|_{u}}{\left|J_{\left(\mu,\nu\right)}\right|}+\Gamma_{u}^{\mathbf{S}}\right)\right]-\frac{\partial\nu}{\partial u}\left[-2\mathbf{m}\cdot\mathbf{n}_{v}-\mathbf{m}\cdot\mathbf{n}\left(-\frac{\left|J_{\left(\mu,\nu\right)}\right|_{v}}{\left|J_{\left(\mu,\nu\right)}\right|}+\Gamma_{v}^{\mathbf{S}}\right)\right]
\end{gather}
and simiarly:
\begin{gather}
\frac{2\left|J_{\left(\mu,\nu\right)}\right|}{\sqrt{\left|\mathbf{M}\right|\left|\mathbf{N}\right|}}\frac{\partial L}{\partial\nu}=-2\mathbf{m}\cdot\left(-\frac{\partial\mu}{\partial v}\mathbf{n}_{u}+\frac{\partial\mu}{\partial u}\mathbf{n}_{v}\right)\qquad\qquad\qquad\qquad\qquad\qquad\qquad\nonumber \\
\qquad\qquad\qquad\qquad\qquad+\mathbf{m}\cdot\mathbf{n}\frac{\partial\mu}{\partial v}\left[-\frac{\left|J_{\left(\mu,\nu\right)}\right|_{u}}{\left|J_{\left(\mu,\nu\right)}\right|}+\Gamma_{u}^{\mathbf{S}}\right]-\mathbf{m}\cdot\mathbf{n}\frac{\partial\mu}{\partial u}\left[-\frac{\left|J_{\left(\mu,\nu\right)}\right|_{v}}{\left|J_{\left(\mu,\nu\right)}\right|}++\Gamma_{v}^{\mathbf{S}}\right]=\nonumber \\
-\frac{\partial\mu}{\partial v}\left[-2\mathbf{m}\cdot\mathbf{n}_{u}-\mathbf{m}\cdot\mathbf{n}\left(-\frac{\left|J_{\left(\mu,\nu\right)}\right|_{u}}{\left|J_{\left(\mu,\nu\right)}\right|}+\Gamma_{u}^{\mathbf{S}}\right)\right]\qquad\qquad\qquad\qquad\qquad\qquad\qquad\\
+\frac{\partial\mu}{\partial u}\left[-2\mathbf{m}\cdot\mathbf{n}_{v}-\mathbf{m}\cdot\mathbf{n}\left(-\frac{\left|J_{\left(\mu,\nu\right)}\right|_{v}}{\left|J_{\left(\mu,\nu\right)}\right|}+\Gamma_{v}^{\mathbf{S}}\right)\right]
\end{gather}

Calculation of the further terms: 
\begin{eqnarray}
-\frac{\partial L}{\partial\mu_{u}} & =\frac{\partial\nu}{\partial v} & \frac{1}{2\sqrt{\left|J_{\left(\mu,\nu\right)}\right|}}\frac{\mathbf{R}_{u}\times\mathbf{R}_{v}}{\sqrt{\left|\mathbf{R}_{u}\times\mathbf{R}_{v}\right|}}\cdot\frac{\mathbf{S}_{\mu}\times\mathbf{S}_{\nu}}{\sqrt{\left|\mathbf{S}_{\mu}\times\mathbf{S}_{\nu}\right|}}\nonumber \\
-\frac{\partial L}{\partial\mu_{v}} & =-\frac{\partial\nu}{\partial u} & \frac{1}{2\sqrt{\left|J_{\left(\mu,\nu\right)}\right|}}\frac{\mathbf{R}_{u}\times\mathbf{R}_{v}}{\sqrt{\left|\mathbf{R}_{u}\times\mathbf{R}_{v}\right|}}\cdot\frac{\mathbf{S}_{\mu}\times\mathbf{S}_{\nu}}{\sqrt{\left|\mathbf{S}_{\mu}\times\mathbf{S}_{\nu}\right|}}\,.\label{eq:partial-result0}
\end{eqnarray}
The explicit dependencies from $\left(\mu,\nu\right)$ can be removed
using (\ref{eq:metric-transformation}):
\begin{eqnarray}
-\frac{\partial L}{\partial\mu_{u}} & =\frac{\partial\nu}{\partial v} & \frac{1}{2\left|J_{\left(\mu,\nu\right)}\right|}\frac{\mathbf{R}_{u}\times\mathbf{R}_{v}}{\sqrt{\left|\mathbf{R}_{u}\times\mathbf{R}_{v}\right|}}\cdot\frac{\mathbf{S}_{u}\times\mathbf{S}_{v}}{\sqrt{\left|\mathbf{S}_{u}\times\mathbf{S}_{v}\right|}}\nonumber \\
-\frac{\partial L}{\partial\mu_{v}} & =-\frac{\partial\nu}{\partial u} & \frac{1}{2\left|J_{\left(\mu,\nu\right)}\right|}\frac{\mathbf{R}_{u}\times\mathbf{R}_{v}}{\sqrt{\left|\mathbf{R}_{u}\times\mathbf{R}_{v}\right|}}\cdot\frac{\mathbf{S}_{\mu}\times\mathbf{S}_{\nu}}{\sqrt{\left|\mathbf{S}_{u}\times\mathbf{S}_{v}\right|}}\,.\label{eq:partial-result1}
\end{eqnarray}
Getting rid of the explicit independency from $\left(\mu,\nu\right)$
in (\ref{eq:partial-result1}) greatly simplify the following calculations
(using (\ref{eq:SRNF_u}) and (\ref{eq:SRNF_v})): 
\begin{gather*}
\frac{2\left|J_{\left(\mu,\nu\right)}\right|}{\sqrt{\left|\mathbf{M}\right|\left|\mathbf{N}\right|}}\left(-\frac{\partial}{\partial u}\frac{\partial L}{\partial\mu_{u}}-\frac{\partial}{\partial\nu}\frac{\partial L}{\partial\mu_{v}}\right)=\qquad\qquad\qquad\qquad\\
\frac{\partial\nu}{\partial v}\left[-\frac{\left|J_{\left(\mu,\nu\right)}\right|_{u}}{\left|J_{\left(\mu,\nu\right)}\right|}\mathbf{m}\cdot\mathbf{n}+\left(\mathbf{m}_{u}+\frac{\mathbf{m}}{2}\Gamma_{u}^{\mathbf{R}}\right)\cdot\mathbf{n}+\left(\mathbf{n}_{u}+\frac{\mathbf{n}}{2}\Gamma_{u}^{\mathbf{S}}\right)\cdot\mathbf{m}\right]\\
-\frac{\partial\nu}{\partial u}\left[-\frac{\left|J_{\left(\mu,\nu\right)}\right|_{v}}{\left|J_{\left(\mu,\nu\right)}\right|}\mathbf{m}\cdot\mathbf{n}+\left(\mathbf{m}_{v}+\frac{\mathbf{m}}{2}\Gamma_{v}^{\mathbf{R}}\right)\cdot\mathbf{n}+\left(\mathbf{n}_{v}+\frac{\mathbf{n}}{2}\Gamma_{v}^{\mathbf{S}}\right)\cdot\mathbf{m}\right]
\end{gather*}
hence the '$\mu$ component' of the Euler-Lagrange equation becomes
\begin{gather*}
\frac{2\left|J_{\left(\mu,\nu\right)}\right|}{\sqrt{\left|\mathbf{M}\right|\left|\mathbf{N}\right|}}\left(\frac{\partial L}{\partial\mu}-\frac{\partial}{\partial u}\frac{\partial L}{\partial\mu_{u}}-\frac{\partial}{\partial\nu}\frac{\partial L}{\partial\mu_{v}}\right)=\qquad\qquad\qquad\qquad\\
\frac{\partial\nu}{\partial v}\left[\mathbf{m}_{u}\cdot\mathbf{n}-\mathbf{n}_{u}\cdot\mathbf{m}+\frac{1}{2}\mathbf{m}\cdot\mathbf{n}\left(\Gamma_{u}^{\mathbf{R}}-\Gamma_{u}^{\mathbf{S}}\right)\right]\\
-\frac{\partial\nu}{\partial u}\left[\mathbf{m}_{v}\cdot\mathbf{n}-\mathbf{n}_{v}\cdot\mathbf{m}+\frac{1}{2}\mathbf{m}\cdot\mathbf{n}\left(\Gamma_{v}^{\mathbf{R}}-\Gamma_{v}^{\mathbf{S}}\right)\right]\,.
\end{gather*}
With similar calculation, the '$\nu$ component' of the Euler-Lagrange
equation is 
\begin{gather*}
\frac{2\left|J_{\left(\mu,\nu\right)}\right|}{\sqrt{\left|\mathbf{M}\right|\left|\mathbf{N}\right|}}\left(\frac{\partial L}{\partial\nu}-\frac{\partial}{\partial u}\frac{\partial L}{\partial\nu_{u}}-\frac{\partial}{\partial\nu}\frac{\partial L}{\partial\nu_{v}}\right)=\qquad\qquad\qquad\qquad\\
-\frac{\partial\mu}{\partial v}\left[\mathbf{m}_{u}\cdot\mathbf{n}-\mathbf{n}_{u}\cdot\mathbf{m}+\frac{1}{2}\mathbf{m}\cdot\mathbf{n}\left(\Gamma_{u}^{\mathbf{R}}-\Gamma_{u}^{\mathbf{S}}\right)\right]\\
\quad\frac{\partial\mu}{\partial u}\left[\mathbf{m}_{v}\cdot\mathbf{n}-\mathbf{n}_{v}\cdot\mathbf{m}+\frac{1}{2}\mathbf{m}\cdot\mathbf{n}\left(\Gamma_{v}^{\mathbf{R}}-\Gamma_{v}^{\mathbf{S}}\right)\right]\,.
\end{gather*}
Summarizing these components, the Euler-Lagrange equation associated
with the reparameterization problem (\ref{eq:srnf-min1}) using the
inverse of the Jacobian of the reparameterization (\ref{eq:inverse-Jacobian})
is
\begin{gather}
2\left[\begin{array}{c}
\frac{\partial L}{\partial\mu}-\frac{\partial}{\partial u}\frac{\partial L}{\partial\mu_{u}}-\frac{\partial}{\partial\nu}\frac{\partial L}{\partial\mu_{v}}\\
\frac{\partial L}{\partial\nu}-\frac{\partial}{\partial u}\frac{\partial L}{\partial\nu_{u}}-\frac{\partial}{\partial\nu}\frac{\partial L}{\partial\nu_{v}}
\end{array}\right]=\qquad\qquad\qquad\qquad\nonumber \\
\qquad\qquad\qquad\qquad\sqrt{\left|\mathbf{M}\right|\left|\mathbf{N}\right|}\left[J_{\left(\mu,\nu\right)}\right]^{-1}\left[\begin{array}{c}
\mathbf{m}_{u}\cdot\mathbf{n}-\mathbf{n}_{u}\cdot\mathbf{m}+\frac{1}{2}\mathbf{m}\cdot\mathbf{n}\left(\Gamma_{u}^{\mathbf{R}}-\Gamma_{u}^{\mathbf{S}}\right)\\
\mathbf{m}_{v}\cdot\mathbf{n}-\mathbf{n}_{v}\cdot\mathbf{m}+\frac{1}{2}\mathbf{m}\cdot\mathbf{n}\left(\Gamma_{v}^{\mathbf{R}}-\Gamma_{v}^{\mathbf{S}}\right)
\end{array}\right]\,.\label{eq:SRNF-Euler-Lagrange}
\end{gather}

\section*{Appendix B\label{sec:Appendix-B}}

The optimal parameterization for the RPSN problem between surfaces
$\mathbf{R}\left(u,v\right)$ (reference) and $\mathbf{S}\left(u,v\right)=\mathbf{S}\left(\mu\left(u,v\right),\nu\left(u,v\right)\right)$
can be formulated as

\begin{equation}
\underset{\mu,\nu}{\min}\iint\left(\mathbf{R}\left(u,v\right)\sqrt{\left|\mathbf{R}_{u}\times\mathbf{R}_{v}\right|}-\mathbf{S}\left(\mu,\nu\right)\sqrt{\left|J_{\left(\mu,\nu\right)}\right|\left|\mathbf{S}_{\mu}\times\mathbf{S}_{\nu}\right|}\right)^{2}dudv\label{eq:rpsn-min1}
\end{equation}
$\left(\left|\mathbf{S}_{u}\times\mathbf{S}_{v}\right|=\left|J_{\left(\mu,\nu\right)}\right|\left|\mathbf{S}_{\mu}\times\mathbf{S}_{\nu}\right|,\:\left|J_{\left(\mu,\nu\right)}\right|=\frac{\partial\mu}{\partial u}\frac{\partial\nu}{\partial v}-\frac{\partial\mu}{\partial v}\frac{\partial\nu}{\partial u}\right)$.
Similarly to the SRNF case, the minimization of (\ref{eq:rpsn-min1})
is equivalent with
\begin{gather}
\underset{\mu,\nu}{\min}-\iint\mathbf{R}\left(u,v\right)\cdot\mathbf{S}\left(\mu,\nu\right)\sqrt{\left|\mathbf{R}_{u}\times\mathbf{R}_{v}\right|}\sqrt{\left|J_{\left(\mu,\nu\right)}\right|\left|\mathbf{S}_{\mu}\times\mathbf{S}_{\nu}\right|}dudv\label{eq:rpsn-min2}
\end{gather}
where the Lagrangian is
\begin{gather}
L=-\sqrt{\left|J_{\left(\mu,\nu\right)}\right|\left|\mathbf{M}\right|\left|\mathbf{S}_{\mu}\times\mathbf{S}_{\nu}\right|}\mathbf{R}\cdot\mathbf{S}\left(\mu,\nu\right)\,.\label{eq:rpsn-Lagrangian}
\end{gather}
$\left(\mathbf{M}=\mathbf{R}_{u}\times\mathbf{R}_{v},\:\mathbf{m}=\frac{\mathbf{M}}{\left|\mathbf{M}\right|}\right)$.
From the latter 
\begin{gather}
\frac{\partial L}{\partial\mu}=-\sqrt{\left|J_{\left(\mu,\nu\right)}\right|\left|\mathbf{M}\right|}\mathbf{R}\cdot\left[\sqrt{\left|\mathbf{N}\right|}\mathbf{S}_{\mu}+\mathbf{S}\frac{\mathbf{n}}{2\sqrt{\left|\mathbf{N}\right|}}\cdot\left(\mathbf{S}_{\mu\mu}\times\mathbf{S}_{\nu}+\mathbf{S}_{\mu}\times\mathbf{S}_{\mu\nu}\right)\right]\nonumber \\
\qquad\qquad\qquad\qquad\qquad\qquad\qquad\qquad=-\sqrt{\left|\mathbf{M}\right|\left|\mathbf{N}\right|}\left(\mathbf{R}\cdot\mathbf{S}_{\mu}+\frac{1}{2}\Gamma_{\mu}^{\mathbf{S}}\mathbf{R}\cdot\mathbf{S}\right)\label{eq:dL/dmu-1}
\end{gather}
$\left(\mathbf{N}=\mathbf{S}_{u}\times\mathbf{S}_{v},\:\mathbf{n}=\frac{\mathbf{N}}{\left|\mathbf{N}\right|}\right)$
and similarly 
\begin{equation}
\frac{\partial L}{\partial\nu}=-\sqrt{\left|\mathbf{M}\right|\left|\mathbf{N}\right|}\left(\mathbf{R}\cdot\mathbf{S}_{\nu}+\frac{1}{2}\Gamma_{\nu}^{\mathbf{S}}\mathbf{R}\cdot\mathbf{S}\right)\label{eq:dL/dnu-1}
\end{equation}
where Christoffel divergences $\Gamma_{\mu}^{\mathbf{S}}$, $\Gamma_{\nu}^{\mathbf{S}}$
are defined in (\ref{eq:Christoffel-divergences}). Using (\ref{eq:inverse-Jacobian})
(\ref{eq:Christoffel-divergence-transform-mu}), equations (\ref{eq:dL/dmu-1})
expressed with the original parameters $\left(u,v\right)$ can be
rearranged as:
\begin{gather}
\frac{2\left|J_{\left(\mu,\nu\right)}\right|}{\sqrt{\left|\mathbf{M}\right|\left|\mathbf{N}\right|}}\frac{\partial L}{\partial\mu}=-2\mathbf{R}\cdot\left(\frac{\partial\nu}{\partial v}\mathbf{S}_{u}-\frac{\partial\nu}{\partial u}\mathbf{S}_{v}\right)\qquad\qquad\qquad\qquad\qquad\qquad\qquad\nonumber \\
\qquad\qquad\qquad\qquad\qquad-\mathbf{R}\cdot\mathbf{S}\frac{\partial\nu}{\partial v}\left[-\frac{\left|J_{\left(\mu,\nu\right)}\right|_{u}}{\left|J_{\left(\mu,\nu\right)}\right|}+\Gamma_{u}^{\mathbf{S}}\right]+\mathbf{R}\cdot\mathbf{S}\frac{\partial\nu}{\partial v}\left[-\frac{\left|J_{\left(\mu,\nu\right)}\right|_{v}}{\left|J_{\left(\mu,\nu\right)}\right|}+\Gamma_{v}^{\mathbf{S}}\right]=\nonumber \\
\frac{\partial\nu}{\partial v}\left[-2\mathbf{R}\cdot\mathbf{S}_{u}-\mathbf{R}\cdot\mathbf{S}\left(-\frac{\left|J_{\left(\mu,\nu\right)}\right|_{u}}{\left|J_{\left(\mu,\nu\right)}\right|}+\Gamma_{u}^{\mathbf{S}}\right)\right]\qquad\qquad\qquad\qquad\qquad\qquad\qquad\\
-\frac{\partial\nu}{\partial u}\left[-2\mathbf{R}\cdot\mathbf{S}_{v}-\mathbf{R}\cdot\mathbf{S}\left(-\frac{\left|J_{\left(\mu,\nu\right)}\right|_{v}}{\left|J_{\left(\mu,\nu\right)}\right|}+\Gamma_{v}^{\mathbf{S}}\right)\right]
\end{gather}
Calculation of the further terms with the removal of the explicit
$\left(\mu,\nu\right)$ dependencies: 
\begin{eqnarray}
-\frac{\partial L}{\partial\mu_{u}} & =\frac{\partial\nu}{\partial v} & \frac{1}{2\left|J_{\left(\mu,\nu\right)}\right|}\sqrt{\left|\mathbf{M}\right|\left|\mathbf{N}\right|}\mathbf{R}\cdot\mathbf{S}\nonumber \\
-\frac{\partial L}{\partial\mu_{v}} & =-\frac{\partial\nu}{\partial u} & \frac{1}{2\left|J_{\left(\mu,\nu\right)}\right|}\sqrt{\left|\mathbf{M}\right|\left|\mathbf{N}\right|}\mathbf{R}\cdot\mathbf{S},
\end{eqnarray}
Then using (\ref{eq:RPSN_u}) and (\ref{eq:RPSN_v})): 
\begin{gather*}
\frac{2\left|J_{\left(\mu,\nu\right)}\right|}{\sqrt{\left|\mathbf{M}\right|\left|\mathbf{N}\right|}}\left(-\frac{\partial}{\partial u}\frac{\partial L}{\partial\mu_{u}}-\frac{\partial}{\partial\nu}\frac{\partial L}{\partial\mu_{v}}\right)=\qquad\qquad\qquad\qquad\\
\frac{\partial\nu}{\partial v}\left[-\frac{\left|J_{\left(\mu,\nu\right)}\right|_{u}}{\left|J_{\left(\mu,\nu\right)}\right|}\mathbf{R}\cdot\mathbf{S}+\left(\mathbf{R}_{u}+\frac{\mathbf{R}}{2}\Gamma_{u}^{\mathbf{R}}\right)\cdot\mathbf{S}+\left(\mathbf{S}_{u}+\frac{\mathbf{S}}{2}\Gamma_{u}^{\mathbf{S}}\right)\cdot\mathbf{R}\right]\\
-\frac{\partial\nu}{\partial u}\left[-\frac{\left|J_{\left(\mu,\nu\right)}\right|_{v}}{\left|J_{\left(\mu,\nu\right)}\right|}\mathbf{R}\cdot\mathbf{S}+\left(\mathbf{R}_{v}+\frac{\mathbf{R}}{2}\Gamma_{v}^{\mathbf{R}}\right)\cdot\mathbf{S}+\left(\mathbf{S}_{v}+\frac{\mathbf{S}}{2}\Gamma_{v}^{\mathbf{S}}\right)\cdot\mathbf{R}\right]
\end{gather*}
hence the '$\mu$ component' of the Euler-Lagrange equation becomes
\begin{gather*}
\frac{2\left|J_{\left(\mu,\nu\right)}\right|}{\sqrt{\left|\mathbf{M}\right|\left|\mathbf{N}\right|}}\left(\frac{\partial L}{\partial\mu}-\frac{\partial}{\partial u}\frac{\partial L}{\partial\mu_{u}}-\frac{\partial}{\partial\nu}\frac{\partial L}{\partial\mu_{v}}\right)=\qquad\qquad\qquad\qquad\\
\frac{\partial\nu}{\partial v}\left[\mathbf{R}_{u}\cdot\mathbf{S}-\mathbf{S}_{u}\cdot\mathbf{R}+\frac{1}{2}\mathbf{R}\cdot\mathbf{S}\left(\Gamma_{u}^{\mathbf{R}}-\Gamma_{u}^{\mathbf{S}}\right)\right]\\
-\frac{\partial\nu}{\partial u}\left[\mathbf{R}_{v}\cdot\mathbf{S}-\mathbf{S}_{v}\cdot\mathbf{R}+\frac{1}{2}\mathbf{R}\cdot\mathbf{S}\left(\Gamma_{v}^{\mathbf{R}}-\Gamma_{v}^{\mathbf{S}}\right)\right]\,.
\end{gather*}
With similar calculation, the '$\nu$ component' of the Euler-Lagrange
equation is 
\begin{gather*}
\frac{2\left|J_{\left(\mu,\nu\right)}\right|}{\sqrt{\left|\mathbf{M}\right|\left|\mathbf{N}\right|}}\left(\frac{\partial L}{\partial\nu}-\frac{\partial}{\partial u}\frac{\partial L}{\partial\nu_{u}}-\frac{\partial}{\partial\nu}\frac{\partial L}{\partial\nu_{v}}\right)=\qquad\qquad\qquad\qquad\\
-\frac{\partial\mu}{\partial v}\left[\mathbf{R}_{u}\cdot\mathbf{S}-\mathbf{S}_{u}\cdot\mathbf{R}+\frac{1}{2}\mathbf{R}\cdot\mathbf{S}\left(\Gamma_{u}^{\mathbf{R}}-\Gamma_{u}^{\mathbf{S}}\right)\right]\\
\quad\frac{\partial\mu}{\partial u}\left[\mathbf{R}_{v}\cdot\mathbf{S}-\mathbf{S}_{v}\cdot\mathbf{R}+\frac{1}{2}\mathbf{R}\cdot\mathbf{S}\left(\Gamma_{v}^{\mathbf{R}}-\Gamma_{v}^{\mathbf{S}}\right)\right]\,.
\end{gather*}
Summarizing these components, the Euler-Lagrange equation associated
with the reparameterization problem (\ref{eq:srnf-min1}) using the
inverse of the Jacobian of the reparameterization (\ref{eq:inverse-Jacobian})
is 
\begin{gather}
2\left[\begin{array}{c}
\frac{\partial L}{\partial\mu}-\frac{\partial}{\partial u}\frac{\partial L}{\partial\mu_{u}}-\frac{\partial}{\partial\nu}\frac{\partial L}{\partial\mu_{v}}\\
\frac{\partial L}{\partial\nu}-\frac{\partial}{\partial u}\frac{\partial L}{\partial\nu_{u}}-\frac{\partial}{\partial\nu}\frac{\partial L}{\partial\nu_{v}}
\end{array}\right]=\qquad\qquad\qquad\qquad\nonumber \\
\qquad\qquad\qquad\qquad\sqrt{\left|\mathbf{M}\right|\left|\mathbf{N}\right|}\left[J_{\left(\mu,\nu\right)}\right]^{-1}\left[\begin{array}{c}
\mathbf{R}_{u}\cdot\mathbf{S}-\mathbf{S}_{u}\cdot\mathbf{R}+\frac{1}{2}\mathbf{R}\cdot\mathbf{S}\left(\Gamma_{u}^{\mathbf{R}}-\Gamma_{u}^{\mathbf{S}}\right)\\
\mathbf{R}_{v}\cdot\mathbf{S}-\mathbf{S}_{v}\cdot\mathbf{R}+\frac{1}{2}\mathbf{R}\cdot\mathbf{S}\left(\Gamma_{v}^{\mathbf{R}}-\Gamma_{v}^{\mathbf{S}}\right)
\end{array}\right]\,.\label{eq:RPSN-Euler-Lagrange}
\end{gather}

\bibliographystyle{plain}
\bibliography{refs}

\end{document}